\documentclass[reqno,12pt]{amsart}
\usepackage{amsmath,amsthm,amscd,amssymb}
\newcommand{\bbC}{{\mathbb{C}}}
\newcommand{\bbD}{{\mathbb{D}}}

\newcommand{\bbG}{{\mathbb{G}}}

\newcommand{\bbR}{{\mathbb{R}}}

\newcommand{\lb}{\label}
\newcommand{\f}{\frac}

\newcommand{\ol}{\overline}

\newcommand{\wti}{\widetilde  }

\newcommand{\spec}{\text{\rm{spec}}}

\newcommand{\s}{\text{\rm{s}}}

\newcommand{\supp}{\text{\rm{supp}}}

\newcommand{\bi}{\bibitem}

\newcommand{\beq}{\begin{equation}}
\newcommand{\eeq}{\end{equation}}
\newcommand{\ba}{\begin{align}}
\newcommand{\ea}{\end{align}}
\newcommand{\veps}{\varepsilon}

\DeclareMathOperator{\Ima}{Im}

\let\det=\undefined\DeclareMathOperator{\det}{det}
\allowdisplaybreaks
\numberwithin{equation}{section}

\newtheorem{theorem}{Theorem}[section]

\theoremstyle{remark}
\newtheorem*{remark}{Remark}

\theoremstyle{definition}
\newtheorem*{definition}{Definition}

\newcommand{\abs}[1]{\lvert#1\rvert}

\newcounter{smalllist}
\newenvironment{SL}{\begin{list}{{\rm\roman{smalllist})}}{%
\setlength{\topsep}{0mm}\setlength{\parsep}{0mm}\setlength{\itemsep}{0mm}%
\setlength{\labelwidth}{2em}\setlength{\leftmargin}{2em}\usecounter{smalllist}%
}}{\end{list}}

\begin{document}

\title[OPs With Exponential Decay]{Orthogonal Polynomials With Exponentially Decaying
Recursion Coefficients}

\author[B.~Simon]{Barry Simon*}

\thanks{$^*$ Mathematics 253-37, California Institute of Technology, Pasadena, CA 91125, USA.
E-mail: bsimon@caltech.edu. Supported in part by NSF grant DMS-0140592}
\thanks{Submitted to the Proceedings of S.~Molchanov's 65th Birthday Conference}

\dedicatory{Dedicated to S.~Molchanov on his 65th birthday}

\date{March 2, 2006}

\begin{abstract} We review recent results on necessary and sufficient conditions for
measures on $\mathbb{R}$ and $\partial\mathbb{D}$ to yield exponential decay of the recursion
coefficients of the corresponding orthogonal polynomials. We include results on the relation
of detailed asymptotics of the recursion coefficients to detailed analyticity of the
measures. We present an analog of Carmona's formula for OPRL. A major role is played
by the Szeg\H{o} and Jost functions.
\end{abstract}

\maketitle

\section{Introduction: Szeg\H{o} and Jost Functions} \lb{s1}

In broad strokes, spectral theory concerns the connection between the coefficients in
differential or difference equations and the spectral measures associated to those
equations. The process of going from coefficients to the measures is the direct
problem, and the other direction is the inverse spectral problem. The gems of spectral
theory are ones that set up one-one correspondences between classes of measures and
coefficients with some properties. Examples are Verblunsky's form of Szeg\H{o}'s
theorem \cite{V36} and the Killip-Simon theorem \cite{KS}. In this paper, our goal is to
describe (mainly) recent results involving such gems for orthogonal polynomials
whose recursion coefficients decay exponentially. These are technically simpler
systems than the $L^2$ results just quoted but have more involved details.

The two classes we discuss are orthogonal polynomials on the real line (OPRL)
and on the unit circle (OPUC). For the OPRL case, we have a probability measure,
$d\rho$, on $\bbR$ of bounded but infinite support whose orthonormal polynomials,
$p_n(x)$, obey
\begin{equation} \lb{1.1}
xp_n(x) = a_{n+1} p_{n+1}(x) + b_{n+1} p_n(x) + a_n p_{n-1}(x)
\end{equation}
with $b_n\in\bbR$ and $a_n\in (0,\infty)$ and called Jacobi parameters.
$\{a_n, b_n\}_{n=1}^\infty$ is a description of $d\rho$ in that there is a one-one
correspondence between bounded sets of such Jacobi parameters and such $d\rho$'s.
For background discussion of OPRL, see \cite{Chi,FrB,OPUC1,Szb}.

For the OPUC case, $\bbD=\{z\in\bbC \mid\abs{z}<1\}$, and $d\mu$ is a probability
measure on $\partial\bbD$ whose support is not a finite set. The orthonormal
polynomials, $\varphi_n(z)$, obey the Szeg\H{o} recursion relation:
\begin{align}
z\varphi_n(z) &= \rho_n \varphi_{n+1}(z) + \bar\alpha_n \varphi_n^*(z) \lb{1.2} \\
\varphi_n^*(z) &= z^n \, \ol{\varphi_n (1/\bar z)} \lb{1.3} \\
\rho_n &= (1-\abs{\alpha_n}^2)^{1/2} \lb{1.4}
\end{align}
with $\alpha_n\in\bbD$ and called Verblunsky coefficients. $\{\alpha_n\}_{n=0}^\infty$
is a description of $d\mu$ in that there is a one-one correspondence between
sequences of $\alpha_n$ obeying $\abs{\alpha_n}< 1$ and such $d\mu$'s. For background
discussion of OPUC, see \cite{GBk,1Foot,OPUC1,OPUC2,Szb}.

The measure theoretic side of the equivalences will be in terms of a derived object,
rather than the measures themselves. For OPUC, the object is $D(z)$, the Szeg\H{o}
function \cite[Section~2.4]{OPUC1}. One says the Szeg\H{o} condition holds if and
only if
\begin{equation} \lb{1.5}
d\mu(\theta) = w(\theta)\, \f{d\theta}{2\pi} + d\mu_\s
\end{equation}
where $d\mu_\s$ is singular and
\begin{equation} \lb{1.6}
\int \log (w(\theta)) \, \f{d\theta}{2\pi} >-\infty
\end{equation}
(which is known to be equivalent to $\sum_{n=0}^\infty \abs{\alpha_n}^2 <\infty$).
In that case, $D(z)$ is defined by
\begin{equation} \lb{1.7}
D(z) =\exp \biggl( \int \f{e^{i\theta}+z}{e^{i\theta}-z}\, \log(w(\theta))\,
\f{d\theta}{4\pi}\biggr)
\end{equation}
which obeys
\begin{equation} \lb{1.8}
\varphi_n^*(z) \to D(z)^{-1}
\end{equation}
if $\abs{z}<1$.

$D(z)$ does not uniquely determine $d\mu$, but it does if $d\mu_\s =0$, as it will
be in our cases of interest, since
\begin{equation} \lb{1.9}
w(\theta) = \lim_{r\uparrow 1}\, \abs{D(re^{i\theta})}^2
\end{equation}
for a.e.\ $\theta$.

For OPRL, the object is the Jost function. The situation is not as clean as the OPUC
case in that there are not simple necessary and sufficient conditions for existence
in terms of the measure. There are necessary and sufficient conditions in terms of the
Jacobi parameters (see \cite{Jost1} and \cite[Section~13.9]{OPUC2}) but not for the measure.
However, there are sufficient conditions for existence that suffice for us here. Suppose
\begin{equation} \lb{1.10}
d\rho(x) = f(x)\, dx + d\rho_\s (x)
\end{equation}
where $f$ is supported on $[-2,2]$, and outside this set, the singular part,
$d\rho_\s$, has only pure points $\{E_j^\pm\}_{j=0}^{N_\pm}$ with
\begin{equation} \lb{1.11}
E_1^- < E_2^- < \cdots < -2 < 2 < \cdots < E_2^+ < E_1^+
\end{equation}
and suppose
\begin{equation} \lb{1.12}
\sum_{j,\pm}\, (\abs{E_j^\pm}-2)^{1/2} < \infty
\end{equation}
and that
\begin{equation} \lb{1.13}
\int_{-2}^2 \log (f(x)) (4-x^2)^{-1/2} \, dx > -\infty
\end{equation}

Then (originally in Peherstorfer--Yuditskii \cite{PY}; see also Simon--Zlato\v{s} \cite{SZ}
and \cite[Theorem~13.8.9]{OPUC2}) there is an analytic function $u$ on $\bbD$ so that
its zeros are precisely those points $z_j^\pm$ in $\bbD$ given by
\begin{equation} \lb{1.14}
z_j^\pm + (z_j^\pm)^{-1} = E_j^\pm
\end{equation}
and if $B$ is the Blaschke product (convergent by \eqref{1.12})
\begin{equation} \lb{1.15}
B(z) = \prod_{j=\pm 1} \, \f{(z-z_j^\pm)}{1-\bar z_j^\pm z}
\end{equation}
then $B u^{-1}\in H^2$ and the boundary values of $u$ obey
\begin{equation} \lb{1.16}
\abs{u(e^{i\theta})}^2 \Ima M(e^{i\theta}) = \sin\theta
\end{equation}
where
\begin{equation} \lb{1.17}
M(z) = \int \f{d\rho(x)}{z+z^{-1} -x}
\end{equation}
(so $\Ima M(e^{i\theta})$ is related to $f(2\cos\theta)$).

These properties determine $u$ uniquely. Unlike the OPUC case, $u$ does not
determine $d\rho$ even if $d\rho_\s \restriction [-2,2]=0$ for $u$ only determines
$f$ and the localization of the pure points of $d\rho$ on $\bbR\backslash [-2,2]$.
To recover $d\rho$, we also need to know the weights of the pure points;
equivalently, the residues of the poles of $M$ at the $z_j^\pm$.

The theme of this review is that detailed results on exponential decay of recursion
coefficients are equivalent to analyticity results on $D^{-1}$ in the OPUC case and
$u$ in the OPRL case. That exponential decay implies analyticity has been in
the physics literature for Schr\"odinger operators for over fifty years. The subtle
aspect is the strict equivalence --- an idea that appeared first in Nevai--Totik
\cite{NT89}.

In Section~\ref{s2}, we discuss some aspects of finite range potentials, and
in Section~\ref{s3}, following Nevai--Totik \cite{NT89} and Damanik--Simon
\cite{DaSim2}, the initial equivalence. In Section~\ref{s4}, following Simon
\cite{S303,S305}, we discuss detailed exponential asymptotics and meromorphic $S$ and $u$.

\smallskip
I would like to thank J.~Christiansen, L.~Golinskii, P.~Nevai, and
V.~Totik for useful comments.

\smallskip
Stas Molchanov is a leading figure in spectral theory. It is a pleasure to
present this birthday bouquet to him.

\section{Finite Range} \lb{s2}

In this section, we present new results on approximation by finite range
``potentials." We begin with an OPRL analog of Carmona's  result \cite{Carm}
on boundary condition averaging for Schr\"odinger operators. We will also see
that Bernstein--Szeg\H{o} measures for OPUC can be viewed through the Carmona
lens. Carmona's proof relies on computing derivatives of Pr\"ufer variables
--- our proof here is spectral averaging making the relation to \cite{SimWolff}
transparent.

Let $J$ be the semi-infinite Jacobi matrix
\begin{equation} \lb{2.1}
J=\begin{pmatrix}
b_1 & a_1 & 0 & \dots \\
a_1 & b_2 & a_2 & \dots \\
0 & a_2 & b_3 & \dots \\
\hdotsfor{4} \\
\hdotsfor{4}
\end{pmatrix}
\end{equation}
associated to an OPRL with measure $d\rho$; $J_{n;F}$, the matrix obtained from
the top $n$ rows and $n$ left columns; and $J_n^b$, the matrix $J_{n;F}$ with
$b_n$ replaced by $b_n +b$. Here $b\in\bbC$. Notice that if $\Ima b\leq 0$, then
$\spec(J_n^b)\subset\bbC\backslash\bbC_+$, so
\[
m_n^{(b)}(z) = \langle\delta_0, (J_n^{(b)}-z)^{-1}\delta_0\rangle
\]
is analytic for $b\in\ol{\bbC}_-$ and $z$ fixed in $\bbC_+$.

If $d\rho$ is a determinate moment problem, then $J$ is essentially selfadjoint on
finite sequences \cite{S270}, so
\begin{equation} \lb{2.2}
m_n^{(b)}(z)\to m(z) \equiv \int \f{d\rho(x)}{x-z}
\end{equation}
for any $b$. Thus, if $d\nu_n$ is defined by
\begin{align}
\wti m_n(z) &\equiv \f{1}{\pi} \int_{-\infty}^\infty m_n^{(b)}(z) \, \f{db}{1+b^2} \lb{2.3} \\
&= \int \f{d\nu^{(n)}(x)}{x-z} \lb{2.4}
\end{align}
then
\begin{equation} \lb{2.5}
d\nu^{(n)}\to d\rho
\end{equation}
weakly. $d\nu^{(n)}$ is thus the average over $b$ of the pure point spectral measures
of $J_n^{(b)}$.

\begin{theorem}\lb{T2.1} If $p_n(x)$ are the orthonormal OPRL, then
\begin{equation} \lb{2.6x}
d\nu^{(n)}(x) = \f{dx}{\pi (a_n^2 p_n^2(x) + p_{n-1}^2(x))}
\end{equation}
\end{theorem}

In particular, the right-hand side of \eqref{2.6x} converges weakly to $d\rho$.
More is true, for Gaussian quadrature implies that if $m_n^{(b)}(z)=\int d\rho_n^{(b)}
(x)(x-z)^{-1}$, then $\int x^\ell \, d\rho_n^{(b)}(x) =\int x^\ell\, d\rho(x)$ for
$\ell\leq 2n-2$, and thus,
\begin{equation} \lb{2.7x}
\int x^\ell\, d\nu^{(n)}(x) = \int x^\ell\, d\rho \qquad
\ell=0, \dots, 2n-2
\end{equation}
Of course, $d\nu^{(n)}$ does not have {\it all\/} moments finite; indeed, $\int
\abs{x}^\ell\, d\nu^{(n)}=\infty$ for $\ell\geq 2n-1$.

\begin{proof}[Proof of Theorem~\ref{T2.1}] It is well known (see \cite[Section~1.2]{OPUC1})
that
\begin{equation} \lb{2.6}
\det (z-J_{n;F}) = P_n(z)
\end{equation}
the monic OPRL, and if $J_{n;F}^{(1)}$ is the matrix obtained by removing the top row
and leftmost column (i.e., $11$ minor), then
\begin{equation} \lb{2.7}
\det (z-J_{n;F}^{(1)}) =Q_n(z)
\end{equation}
the monic second kind polynomial of degree $n-1$.

By expanding $\det(z-J_n^{(b)})$ in minors, we see
\begin{align}
\det(z-J_n^{(b)}) &= P_n(z) -bP_{n-1}(z)  \lb{2.8} \\
&= (a_1\dots a_{n-1}) (a_n p_n(z) - bp_{n-1}(z)) \lb{2.9}
\end{align}
and thus,
\begin{equation} \lb{2.10}
m_n^{(b)}(z) = -\f{(a_n q_n(z) -bq_{n-1}(z))}{(a_n p_n (z) - bp_{n-1}(z))}
\end{equation}

As noted above, if $\Ima b\leq 0$, $m_n^{(b)}(z)$ has its poles in $\Ima z\leq 0$
and thus, if $\Ima z>0$, $m_n^{(b)}(z)$ is analytic in $\Ima b\leq 0$. Thus, we can
close the contour in the lower half-plane and find for $\Ima z> 0$,
\begin{align}
\wti m_n(z) &= m_n^{(b=-i)}(z) \notag \\
&= - \f{(a_n q_n(z) + iq_{n-1}(b))}{(a_n p_n(z) + ip_{n-1}(z))} \lb{2.11}
\end{align}
Thus, $\wti m_n$ is analytic on $\ol{\bbC}_+$, so
\[
d\nu_n(x) = \pi^{-1} \Ima \wti m_n(x)\, dx
\]

Since $p_n$, $p_{n-1}$, $q_n$, $q_{n-1}$ are real on $\bbR$,
\begin{equation} \lb{2.12}
\Ima \wti m_n(x) = \f{a_n (p_{n-1} (x) q_n(x) - p_n(x) q_{n-1}(x))}
{(a_n^2 p_n(x)^2 + p_{n-1}(x)^2)}
\end{equation}
which is \eqref{2.6} by a standard Wronskian calculation (see (1.2.51) of \cite{OPUC1}).
\end{proof}

By this same calculation, one can recover Carmona's formula for the Schr\"odinger
operator case.

One can ask about the analog of this for OPUC. Given a nontrivial measure, $d\mu$, on
$\partial\bbD$ and $\omega =e^{i\theta} \in\partial\bbD$, we define $d\mu_n^{(\omega)}$
to be the trivial measure with Verblunsky coefficients
\begin{align*}
\alpha_j &= \alpha_j (d\mu) \qquad j=0, \dots, n-1 \\
\alpha_n &=\omega
\end{align*}
Then $d\mu_n^{(\omega)}$ is the measure with $n+1$ pure points at the zeros of the
paraorthogonal polynomial (POPUC),
\begin{equation} \lb{2.13}
\Phi_{n+1}^{(\omega)}(z) = z\Phi_n (z) -\bar\omega \Phi_n^*(z)
\end{equation}

\begin{theorem}\lb{T2.2} $d\mu_n\equiv \int \f{d\theta}{2\pi}\, d\mu_n (e^{i\theta})$
is the Bernstein--Szeg\H{o} measure
\begin{equation} \lb{2.14}
d\mu_n = \f{d\theta}{2\pi \abs{\varphi_n (e^{i\theta})}^2}
\end{equation}
\end{theorem}

\begin{proof} If $\psi_n$ are the second kind polynomials, Geronimus' formula for $F(z)$
(see \cite[Theorem~3.2.4]{OPUC1}) implies ($F(z;d\mu)=\int \f{e^{i\theta}+z}{e^{i\theta}-z}
d\mu(\theta)$)
\begin{equation} \lb{2.15}
F(z;d\mu_n^{(\omega)}) = \f{\psi_n^* (z) - \omega z\psi_n(z)}
{\varphi_n^*(z) - \omega z\varphi_n(z)}
\end{equation}
Averaging $\omega$ over $\f{d\theta}{2\pi}$ gives the value at $\omega =0$ since this
function is analytic in $\omega$ for $z$ fixed in $\bbD$. It follows that
\begin{equation} \lb{2.16x}
F(z;d\mu_n)=\f{\psi_n^*(z)}{\varphi_n^*(z)}
\end{equation}
and yields \eqref{2.14} by (3.2.35) of \cite{OPUC1}.
\end{proof}

The Bernstein--Szeg\H{o} approximation also has the property of being the measure
associated to extending the $\alpha$'s up to $n$ to be free beyond $n$ (i.e.,
$\alpha_j=0$ for $j\geq n$). One can ask about the analogous approximation for OPRL. We
will get the function $S_n$ used by Dombrowski--Nevai \cite{DN}:

Let $J_\ell$ be the Jacobi matrix with parameters
\begin{align}
a_n (J_\ell) &= \begin{cases}
a_n(J) & n=1, \dots, \ell-1 \\
1 & n\geq \ell
\end{cases} \lb{2.16} \\
b_n(J_\ell) &= \begin{cases} b_n(J) & n=1, \dots, \ell \\
0 & n\geq \ell \end{cases} \lb{2.17}
\end{align}
According to Theorem~13.6.1 (with $a_\ell$ replaced by $1$), its Jost function is
($x=z+ 1/z$)
\begin{equation} \lb{2.18}
g_\ell (z) = z^\ell \biggl(p_\ell \biggl( z+\f{1}{z}\biggr) -
z p_{\ell-1} \biggl( z + \f{1}{z}\biggr)\biggr)
\end{equation}

Define $S_\ell(x)$ by
\begin{equation} \lb{2.19}
S_\ell \biggl( z+\f{1}{z}\biggr) \equiv g_\ell(z) g_\ell \biggl(\f{1}{z}\biggr)
\end{equation}
Then, by \eqref{2.18},
\begin{equation} \lb{2.20}
S_\ell(x) = p_\ell(x)^2 + p_{\ell-1}(x)^2 - x p_\ell (x) p_{\ell-1}(x)
\end{equation}
Taking into account the different normalization (for us, ``free" is $a_k=1$;
for them, $a_k=\f12$), this is the function $S_\ell(x)$ of Dombrowski--Nevai
\cite{DN}. The approximating measure has a.c.\ part related to $dx/\abs{g_\ell(z)}^2$
on $[-2,2]$ which is $dx/S_\ell(x)$. The eigenvalues of $J_\ell$ are zeros of $S_\ell(x)$
but not all zeros since $S_\ell$ also vanishes if $g_\ell (1/z)=0$, that is,
at antibound state and resonance energies.

For most purposes, \eqref{2.6} is a more useful representation than the one associated
to $S_\ell$.

\section{Necessary and Sufficient Conditions on Exponential Decay} \lb{s3}

The starting point of the recent results on exponential decay is the following
result of Nevai--Totik  for OPUC:

\begin{theorem}[\cite{NT89}] \lb{T3.1} Let $d\mu$ be a nontrivial probability
measure on $\partial\bbD$ and $R>1$. Then the following are equivalent:
\begin{SL}
\item[{\rm{(a)}}]
\begin{equation} \lb{3.1}
\limsup_{n\to\infty} \, \abs{\alpha_n (d\mu)}^{1/n} \leq R^{-1}
\end{equation}
\item[{\rm{(b)}}] $d\mu_\s =0$, the Szeg\H{o} condition \eqref{1.6} holds, and
$D(z)^{-1}$ has an analytic continuation to $\{z\mid \abs{z}<R\}$.
\end{SL}
\end{theorem}

\begin{remark} Since $R^{-1}<1$, \eqref{3.1} is an expression of exponential decay.
\end{remark}

The proof is easy. If \eqref{3.1} holds, Szeg\H{o} recursion first implies inductively
that for $\abs{z} =1$,
\begin{equation} \lb{3.2}
\abs{\Phi_{n+1}(e^{i\theta})}\leq (1+\abs{\alpha_n}) \abs{\Phi_n (e^{i\theta})}
\end{equation}
so
\begin{align}
\sup_{n,\abs{z}\leq 1}\, \abs{\Phi_n^*(z)} &=
\sup_{n,\theta}\, \abs{\Phi_n^*(e^{i\theta})} \quad \text{(by the maximum principle)} \notag \\
&\leq \prod_{j=0}^\infty (1+\abs{\alpha_j}) \equiv C<\infty \lb{3.3}
\end{align}
and thus, for $\abs{z} >1$,
\begin{equation} \lb{3.4}
\abs{\Phi_n(z)} \leq C\abs{z}^n
\end{equation}

Iterating
\begin{equation} \lb{3.5}
\Phi_{n+1}^*(z) = \Phi_n^*(z) -\alpha_n z\Phi_n(z)
\end{equation}
we get
\begin{equation} \lb{3.6}
\Phi_n^*(z) = 1-\sum_{j=0}^{n-1} \alpha_j z\Phi_j(z)
\end{equation}
\eqref{3.1}, \eqref{3.4}, and \eqref{3.6} imply that for any $\veps >0$,
\[
\sup_{n,\abs{z}<R-\veps}\, \abs{\Phi_n^*(z)} <\infty
\]
which implies that $\varphi_n^*(z)$ has a limit for $\abs{z}<R$. This limit defines
the analytic continuation of $D(z)^{-1}$.

For the other direction, one can use either of two similar-looking but distinct
formulae relating $D$ to $\alpha_n$. One can use a formula of Geronimus \cite{GBk}
and Freud \cite{FrB} as Nevai--Totik \cite{NT89} do (it requires $d\mu_\s =0)$
\begin{equation} \lb{3.7}
\alpha_n = -\kappa_\infty \int \ol{\Phi_{n+1}(e^{i\theta})}\, D(e^{i\theta})^{-1}
d\mu (\theta)
\end{equation}
or the following formula of Simon \cite{S303} derived from iterated Szeg\H{o} recursion:
\begin{equation} \lb{3.8}
\alpha_n =-\kappa_\infty^{-1} \kappa_n^2 \int \ol{\Phi_n (e^{i\theta})}\,
[D(e^{i\theta})^{-1}- D(0)^{-1}] e^{-i\theta}\, d\mu(\theta)
\end{equation}
In these formulae,
\[
\kappa_n = \prod_{j=0}^{n-1} (1-\abs{\alpha_j}^2)^{-1/2} \qquad
\kappa_\infty = \lim_{n\to\infty} \, \kappa_n
\]

To get exponential decay of $\alpha_n$ from \eqref{3.7} or \eqref{3.8}, one uses
$\int \ol{\Phi_n (e^{i\theta})} e^{-ij\theta}\, d\mu(\theta)=0$ for $j<n$ and the
Taylor series for $D^{-1}$ to see that $\alpha_n$ is bounded by the tail of the
Taylor series of $D(z)^{-1}$ which, of course, decays exponentially if $D(z)^{-1}$
is analytic in $\abs{z}<R$.

For OPRL, the analogs of Theorem~\ref{T3.1} are due to Damanik--Simon \cite{DaSim2}.
The result is simpler if there are no bound states or resonances where

\begin{definition} We say a measure $d\rho$ on $\bbR$ has no bound states or
resonances if
\begin{equation} \lb{3.9}
d\rho(x) = f(x)\, dx + d\rho_\s
\end{equation}
where
\begin{equation} \lb{3.10}
\supp (d\rho) \subset [-2,2]
\end{equation}
and
\begin{equation} \lb{3.11}
\int (4-x^2)^{-1} f(x)\, dx <\infty
\end{equation}
\end{definition}

\begin{theorem}[\cite{DaSim2}]\lb{T3.2} Let $R>1$. Suppose $d\rho$ has no bound states
or resonances. Then $u(z)$ has an analytic continuation to $\{z\mid\abs{z}<R\}$ if and
only if
\begin{equation} \lb{3.12x}
\limsup [\abs{a_n (d\rho)-1} + \abs{b_n (d\rho)}]^{1/2n} \leq R^{-1}
\end{equation}
\end{theorem}

\cite{DaSim2} has several proofs, but the simplest one is in \cite{S305}. When
\eqref{3.11} holds, there is a measure $d\mu$ on $\partial\bbD$ given by
\begin{equation} \lb{3.12}
d\mu(\theta) = w(\theta) \, \f{d\theta}{2\pi} + d\mu_\s
\end{equation}
where
\begin{equation} \lb{3.13}
w\biggl( \arccos \biggl(\f{x}{2}\biggr)\biggr) =  c (4-x^2)^{-1/2} f(x)
\end{equation}
for suitable $c$ and $d\mu_\s$. The Verblunsky coefficients $\alpha_n$ for $d\mu$
and Jacobi parameters for $d\rho$ are related by (\cite{BCG1,KilNen};
\cite[Section~13.2]{OPUC2})
\begin{align}
b_{n+1} &= \alpha_{2n} -\alpha_{2n+2} -\alpha_{2n+1} (\alpha_{2n} + \alpha_{2n+2}) \lb{3.14} \\
a_{n+1}^2 -1 &= \alpha_{2n+1} -\alpha_{2n+3} -\alpha_{2n+2}^2 (1-\alpha_{2n+3})
(1+\alpha_{2n+1}) - \alpha_{2n+3} \alpha_{2n+1} \lb{3.15}
\end{align}
and the Jost function for $d\rho$ and Szeg\H{o} function for $d\mu$ by
\begin{equation} \lb{3.16}
u(z) =(1-\abs{\alpha_0}^2)(1-\alpha_1) D(z)^{-1}
\end{equation}

From this, it is easy to derive Theorem~\ref{T3.2} from Theorem~\ref{T3.1}.

To understand the situation when $J$ has bound states, we note the analytic continuation
of \eqref{1.16} says
\begin{equation} \lb{3.17}
u(z) u\biggl( \f{1}{z}\biggr) \biggl[M(z)  - M\biggl(\f{1}{z}\biggr)\biggr]
=z-z^{-1}
\end{equation}
(this uses also $u,M$ real on $\bbR$). Recall that if $z_0\in\bbD$ is such that
$z_0 +z_0^{-1}$ is an eigenvalue of $J$, then $u(z_0)=0$. An argument shows that if
$\abs{z_0}>R^{-1}$ and $\abs{a_n-1}+\abs{b_n}\leq CR^{-2n}$, then $u(z_0^{-1})\neq 0$ and
$M(1/z)$ is regular at $z_0$. Thus, \eqref{3.17} implies a relation between $u'(z_0)$,
$u(1/z_0)$, and the residue of the pole of $M(z)$ at $z_0$. This leads to

\begin{definition} Suppose $u$ is analytic in $\{z\mid\abs{z}<R\}$ for some $R>1$ and
$z_0\in\bbD$ with $u(z_0)=0$ and $\abs{z_0}>R^{-1}$. We say the weight of the point mass
at $z_0 + z_0^{-1}$ is canonical if
\begin{equation} \lb{3.18}
\lim_{z\to z_0}\, (z-z_0) M(z_0) = (z_0 -z_0^{-1})
\biggl[u'(z_0) u\biggl( \f{1}{z_0}\biggr)\biggr]
\end{equation}
\end{definition}

\begin{theorem}[\cite{DaSim2}]\lb{T3.3} Fix $R>1$. Then \eqref{3.12} holds if and only if
\begin{SL}
\item[{\rm{(i)}}] $u(z)$ has an analytic continuation to $\{z\mid\abs{z}<R\}$.
\item[{\rm{(ii)}}] The point mass at each $z_0\in\bbD$ with $\abs{z_0} >R^{-1}$
and $u(z_0)=0$ is a canonical weight.
\end{SL}
\end{theorem}

If $u$ is entire and has $m$ zeros in $\bbD$, the set of measures with that $u$ has dimension
$m-1$. A single point on this space has decay at rate faster than any exponential.
Similarly, if $u$ is a polynomial, $\{a_n-1,b_n\}$ has finite support if and only if all
weights are canonical.

\section{Detailed Asymptotics} \lb{s4}

Let $S$ be defined by
\begin{equation} \lb{4.1}
S(z)=-\sum_{j=0}^\infty \alpha_{j-1} z^j
\end{equation}
where $\alpha_{-1}=-1$. Of course, when $D$ exists, both $D(z)^{-1}$ and $S(z)$ are analytic
near $z=0$. Theorem~\ref{T3.1} can be rephrased.

\begin{theorem}\lb{T4.1} The Taylor series of $D(z)^{-1}$ and $S(z)$ have the same radius of
convergence.
\end{theorem}

Barrios, L\'opez, and Saff \cite{BLS} extend this to show $S(z)$ is meromorphic in $\{z\mid
\abs{z}<R+\veps\}$ with a single simple pole at $z=R$ if and only if $D(z)^{-1}$ is
meromorphic in a similar region. This condition on $S$ is, of course, equivalent to
\begin{equation} \lb{4.2}
\alpha_n =CR^{-n} + O(R^{-n(1+\delta)})
\end{equation}
which is how they phrased their result. To go further, it is useful to define
\begin{equation} \lb{4.3}
r(z) = \ol{D(1/\bar z)}\, D(z)^{-1}
\end{equation}
which is analytic in $\{z\mid 1-\veps < \abs{z}<R\}$ if \eqref{3.1} holds. Simon \cite{OPUC1}
proved that $r(z)-S(z)$ is analytic in $\{z\mid 1-\veps <\abs{z} <R^2\}$ when \eqref{3.1}
holds, thereby generalizing \cite{BLS}. The ultimate result of this genre was found independently
by Deift--Ostensson \cite{DO} and Mart\'inez-Finkelshtein et al.\ \cite{MFMS}; an alternate
proof was then found by Simon \cite{S303}.

\begin{theorem} \lb{T4.2} If \eqref{3.1} holds for some $R>1$, then $r(z)-S(z)$ is analytic
in $\{z\mid 1-\veps <\abs{z}<R^3\}$.
\end{theorem}

This is optimal in that there are examples \cite{MFMS,S303} where $S$ (and $r$)
have a simple pole at $z=R$ but $S-r$ has a pole at $z=R^3$.

Motivated by this, Simon \cite{S303} proved:

\begin{theorem} \lb{T4.3} $S(z)$ is an entire meromorphic function if and only if $D(z)^{-1}$ is.
\end{theorem}

One can even relate the poles. Given a set $S$ in $\{z\mid \abs{z}>1\}$ which is discrete,
one defines $\bbG(S)$ to be the set of all products $z_1 \dots z_{n+1} \bar z_{n+2}
\dots \bar z_{2n+1}$ where $z_j\in S$. Then

\begin{theorem}[\cite{S303}] \lb{T4.4} Let $S(z)$ be entire meromorphic and let $P$ be
the poles of $D(z)^{-1}$ and $T$ the set of poles of $S(z)$. Then $P\subset \bbG(T)$
and $T\subset \bbG(P)$.
\end{theorem}

Simon \cite{S305} studies analogs of the results for OPRL. In the Jacobi case, define
\begin{equation} \lb{4.4}
B(z) =1-\sum_{n=0}^\infty [b_{n+1} z^{2n+1} + (a_{n+2} -1) z^{2n+2}]
\end{equation}
The analog of Theorem~\ref{T4.2} is

\begin{theorem}[\cite{S305}] \lb{T4.5} Suppose $R>1$ and
\[
\limsup_{n\to\infty}\, (\abs{a_n^2-1} + \abs{b_n})^{1/2n} = R^{-1}
\]
Then $(1-z^2) u(z) + z^2 u(1/z) B(z)$ is analytic in $\{z\mid R^{-1} < \abs{z} < R^2\}$.
\end{theorem}

As explained there, $R^2$ is optimal. The analog of Theorem~\ref{T4.3} is

\begin{theorem}\lb{T4.6} $B(z)$ is an entire meromorphic function if and only if $u(z)$ is.
\end{theorem}

The connection between poles, that is, the analog of Theorem~\ref{T4.4} is complicated
but appears in \cite{S305}.

\bigskip

\end{document}